\documentclass[reqno, 12pt,a4]{amsart}

\usepackage[mathscr]{eucal}
\usepackage{amssymb}
\usepackage{latexsym}
\usepackage{amsthm}
\usepackage{amscd}
\usepackage{amsmath}

\numberwithin{equation}{section}

\theoremstyle{plain}

\newtheorem{theorem}{Theorem}[section]
\newtheorem{conjecture}[theorem]{Conjecture}

\newtheorem{corollary}[theorem]{Corollary}
\newtheorem{lemma}[theorem]{Lemma}

\theoremstyle{definition}

\newcommand{\bx}{\mbox{\boldmath $x$}}

\def\p{p}
\newcommand{\Z}{{\mathbb{Z}}}
\newcommand{\lm}{{\lambda}}
\newcommand{\qint}[1]{{[p]}_{#1}}
\newcommand{\qnum}[1]{({#1})_{[p]}}
\newcommand{\gstep}{g}
\newcommand\KLR{{\widetilde{H}_n(\zeta)}}

\def\ord{{{o}}}
\def\core{{\mathrm{Core}}}
\def\pminus{{\setminus}}
\def\pplus{{\sqcup}}

\newcommand{\floor}[1]{{\left\lfloor{#1}\right\rfloor}}
\newcommand{\mmu}{{\underline{\mu}}}
\newcommand{\mlm}{{\underline{\lm}}}

\newcommand\coren{c}
\newcommand\cart{C_n}

\def\bd{\Delta}
\newcommand{\Pless}{P_{(\p)}(\leq r)}
\def\inv{{}^{-1}}
\newcommand\frev[2]{\frac{#2}{#1}}
\title[Graded Cartan matrices]{Combinatorics for 
 graded Cartan matrices of the Iwahori-Hecke algebra of type $A$}
\author{Masanori Ando, Takeshi Suzuki 
and Hiro-Fumi Yamada}
\address{Department of Mathematics, Okayama University, Okayama 700-8530, 
Japan}
\email{m\_ando@math.okayama-u.ac.jp;  suzuki@math.okayama-u.ac.jp; 
\newline 
\indent 
yamada@math.okayama-u.ac.jp}
\begin{document}
\begin{abstract}
A $q$-analogue of 
combinatorics concerning the Cartan matrix for the Iwahori-Hecke algebra
of type $A$ is investigated.
We give several descriptions for the determinant of the graded Cartan matrix,
which imply some combinatorial identities.
A conjectural expression for the elementary divisors is
also presented.
\end{abstract}
\maketitle
\section{Introduction}
Let 
$H_n(\zeta)$ be the Iwahori-Hecke algebra of 
type $A_{n-1}$ with the parameter $\zeta$ 
being a primitive $p$-th root of unity.
The irreducible representations of $H_n(\zeta)$ are
labeled by the set $P^{(p)}(n)$ of the
$p$-regular partitions.
The square matrix $C_n=([P(\lambda):D(\mu)])_{\lambda,\mu\in P^{(p)}(n)}$
is called the Cartan matrix of $H_n(\zeta)$,
 where $[P(\lambda):D(\mu)]$ denotes the multiplicity of the irreducible module
$D(\mu)$ 
in a composition series of the 
projective cover $P(\lambda)$ of $D(\lambda)$.
As is well-known, the Cartan matrix can be 
expressed as
$C_n={}^tD_nD_n$,
where $D_n$ is the decomposition matrix of $H_n(\zeta)$.

There are several combinatorial expressions for the elementary divisors and the
determinant of $C_n$.  For example,
when $p$ is prime, it is classically known that
the elementary divisors of 
$C_n$ coincide with
those of the Cartan matrix 
of the symmetric group $\mathfrak{S}_n$ at characteristic $p$, 
and they are given by
\begin{equation}
\left\{\left.
\prod_{i\geq 1}(m_i!)_p
\ \right|\ 
\lambda= (1^{m_1}2^{m_2}\dots )
\in P_{(p)}(n)
\right\},
\end{equation}
where
$P_{(p)}(n)$ denotes the set 
of the $p$-class regular partitions of $n$ 
and
$(k)_p$ is the $p$-part of $k$,
i.e., $(k)_p=p^j$ for $k=a p^j$ with 
$p\nmid a$,
and hence 
\begin{equation}\label{eq;det_clasical}
\det C_n=\prod_{\lm\in P_{(p)}(n)}
\prod_{i\geq 1}(m_i!)_p.
\end{equation}
(see e.g. \cite{nt}). 
%
A bijection called the {\it Glaisher correspondence} \\
$\gamma: P_{(p)}(n)\to P^{(p)}(n)$ gives  
an alternative expression of the elementary divisors:
\begin{equation}
\label{eq;E=G_classical}
 \prod_{i\geq 1}(m_i!)_p = 
p^{\frac{\ell(\lambda)-\ell(\gamma(\lambda))}{p-1}}
 \end{equation}
for $\lambda = (1^{m_1} 2^{m_2} \ldots ) \in P_{(p)}(n)$, 
 where $\ell(\lm)$ denotes the length of the partition
\cite{uy}.

The purpose of the present paper is to give a $q$-analogue of 
combinatorics concerning the Cartan matrix for $H_n(\zeta)$,
and to give some combinatorial identities 
for partitions.
In particular we give a $q$-analogue of both sides of (1.3) and show that their products
over $P_{(p)}(n)$ coincide.
It will turn out that these products equal the determinant
of the $q$-Cartan matrix $C_n(q)$ in the sense of Lascoux, Leclerc and Thibon \cite{llt}.

We describe the composition of this article. 
In Section \ref{sec;partitions}, we recall preliminary results on partitions.
In Section \ref{sec;ed}, we define two weights $w_E$ and $w_H$ for partitions.
Here $w_E$ is a natural $q$-analogue of (1.1), while
$w_H$ is motivated by the work of Hill \cite{h} on block elementary divisors.
We will prove that certain products of $w_E$ and $w_H$ coincide (Corollary 3.4). 
In Section \ref{sec;comb}, we define yet another weight $w_G$ for partitions 
which comes from the Glaisher correspondence.
Also we recall an expression of the block determinant of $C_n(q)$ due
to Tsuchioka (unpublished).
The main result of 
this section (Theorem \ref{th;comb2}) reads
$$
\det C_n(q) 
=\prod_{\lm\in P_{(\p)}(n)}w_G(\lm)
=\prod_{\lm
\in P_{(\p)}(n)}w_E(\lm).
$$
Looking at the equality term-by-term, we obtain an interesting 
partition identity (Corollary \ref{cor;comb2}).
We also give alternative description for Tsuchioka's expression
of the block determinant (Theorem \ref{th;comb1}), 
and present another partition identity 
(Corollary \ref{cor;comb1}).
Theorem \ref{th;comb2} and Theorem \ref{th;comb1} 
are proved in Section \ref{sec;glaisher} 
and Section \ref{sec;pfcomb1}, respectively.
In Section \ref{sec;p=2}, we focus on the case $p=2$. In this case block splitting 
is easily described by means of {\it{H-abacus}} (\cite{uy}).

Representation-theoretic meaning of
the $q$-Cartan matrix $C_n(q)$ 
is best understood if we consider the Khovanov-Lauda-Rouquier algebra
$\KLR$ (\cite {kl1, kl2} and \cite {rouq}).
It is shown by 
Brundan and Kleshchev \cite{bk2} 
that 
$C_n(q)$ 
is the corresponding graded Cartan matrix.
Section 8 is devoted to explaining these relationships.
We conjecture that, when $p$ is prime, 
the diagonal matrix with entries $w_H$ 
has the same elementary divisors with 
the graded Cartan matrix (Conjecture \ref{con}).

\newpage
\section{Partitions}\label{sec;partitions}
%
For a partition $\lm$, 
we let  $m_i(\lm)$ denote the multiplicity of $i$ as its part, and 
we represent $\lm$ as 
 $(1^{m_1(\lm)}2^{m_2(\lm)}\dots)$.

Let $p$ be a fixed integer greater than 1.
A partition $\lm=(1^{m_1}2^{m_2}\dots)$
is said to be $p$-{\it regular} if $m_i <p$ for all $i$, and is said to be $p$-{\it class regular} if
$m_i=0$ for $i$ which is divisible by $p$.  
We let $P(n)$, $P^{(p)}(n)$ and $P_{(p)}(n)$
denote the set of the partitions of $n$,
the
$p$-regular partitions of $n$,  and the set of 
the $p$-class regular partitions of $n$, respectively.
It is well-known (see e.g. \cite{and} that $P^{(p)}(n)$ and $P_{(p)}(n)$ have the same cardinality.

A partition is called a $p$-{\it core} if it has no $p$-hooks.
Let $\core_\p(d)$ denote the set of all $\p$-cores in 
$P(d)$, and let $\coren_p(d)=\sharp \core_p(d)$.

The set $P^{(\p)}(n)$ (and $P_{(\p)}(n)$) labels 
the the set of isomorphism classes of
the irreducible representations of the Iwahori-Hecke algebra
${H}_n(\zeta)$
associated with the symmetric group with the parameter
$\zeta$ being a primitive $p$-th root of $1$.
The set $\sqcup_{0\leq d\leq 
\floor{\frac{n}{\p}}}  \core_{\p}(n-\p d)$
labels the blocks of the Cartan matrix $\cart$
of ${H}_n(\zeta)$,
where $\floor{a}$ denotes the largest integer 
which is not greater than $a$.

For each $p$-core in $ \core_{\p}(n-\p d)$,
the size of the corresponding block matrix  
equals the 
cardinality of the set 
\begin{align*}
&M_{\p-1}(d)=
\left\{\underline\lm=(\lm^{(1)},
\dots,
\lm^{(\p-1)})\mid \lm^{(i)}\in P(d_i),\ 
\Sigma_{i=1}^{\p-1} d_i=d\right\}
\end{align*}
of $(\p-1)$-multipartitions of $d$.
Put
$$
Q_\p(n)=
\bigsqcup_{0\leq d\leq 
\floor{\frac{n}{\p}}}
M_{\p-1}(d)\times \core_{\p}(n-\p d).$$
Then
$\sharp P_{(\p)}(n)=\sharp Q_\p(n)$ 
as they both equal the size of the Cartan matrix $C_n$.

Put
$$\phi(x)=\prod_{n\geq 1}
{(1-x^n)}.$$
Then we have 
$\sum_{n\geq0}
\sharp P(n)x^n=\frac{1}{\phi(x)}$
and the following formulas for
generating functions:
\begin{align}
\sum_{n\geq0}
\sharp P_{(\p)}(n)x^n
&=\frac{\phi(x^\p)}{\phi(x)},
\label{eq;gen_Pp}\\
\sum_{n\geq0}\sharp M_{\p}(n) x^n&=
\frac{1}{{\phi(x)^{\p}}},
\label{eq;gen_M}
\\
\sum_{n\geq0}
\coren_p(n) 
x^n&=\frac{\phi(x^\p)^\p}{\phi(x)}.
\label{eq;gen_core}
\end{align}

We end this section with presenting a useful formula of generating functions.
Let $A(n)$ and $B(n)$ be finite sets indexed by non-negative integers $n$,
and $a(n),\,\,b(n)$ their cardinalities, respectively.  For a positive integer $p$,
consider the set
$$Z(n) = \bigcup_{k\geq0} A(k) \times B(n-pk).$$
Then the generating function of $z(n) = \sharp Z(n)$ is given by
\begin{equation}
\sum_{n\geq 0} z(n) x^{n} 
= \left(\sum_{n\geq 0}a(n)x^{pn}\right)\left(\sum_{n\geq 0}b(n)x^{n}\right).
\end{equation}
Note that equation (2.2) can be shown by this formula.
\section{$p$-part and elementary divisors}
\label{sec;ed}
Let $q$ be an indeterminate. 
For a positive integer $l$, we define $[p]_l$ as the $q$-integer of $p$ with $q^{2l}$ base;
$$[p]_l=\frac{1-q^{2pl}}{1-q^{2l}}
=1+q^{2l}+\dots+q^{2l(p-1)}.
$$
If a positive integer $k$ is written as $k=ap^b$ with 
$\p\nmid a$, 
then $(k)_p = p^b$ is called the $p$-part of $k$,
 and $(k)_{p'} = a$ is called the $p'$-part of $k$.
Put 
$$\qnum{k}=\qint{a}\qint{ap} \cdots \qint{ap^{b-1}},$$ 
which might be called
the {\it{graded}} $p$-part of $k$. 

For a positive integer $m$ with $\p$-adic expansion
$$m=c_0+c_1\p+\dots+c_r\p^r\quad(0\leq c_0,\dots,c_r
\leq\p-1,\ c_r\neq0),$$
define
$$\ord_\p (m)=r+1.$$
The following two lemmas are  verified by direct computations.
\begin{lemma}\label{lem;1}
Let $m\in\Z_{\geq1},\ \p \in \Z_{\geq2}$. 
Then
$$\prod_{j=1}^m\qint{j}^{\ord_\p
\left(
\left\lfloor
\frac{m}{j}\right\rfloor\right)}
=\prod_{j=1}^{m}
[\p ]_j (j)_{[\p]}.
%
$$
\end{lemma}
\begin{lemma}\label{lem;2}
Let $m\in\Z_{\geq1},\ \p \in \Z_{\geq2}$. 
Then
$$
\prod_{j=1}^m(j)_{[\p]}
=
\prod_{j=1}^{\left\lfloor\frac{m}{\p}
\right\rfloor}
[\p ]_j (j)_{[\p]}.
%
$$
\end{lemma}
For a partition $\lm$, we put
\begin{align*}
w_E(\lm)&=\prod_{i\geq 1;\,p\nmid i}\,\,\prod_{j=1}^{m_i(\lm)}\,\qnum{j},\\
w_H(\lm)&=\prod_{j\geq1} \prod_{i\geq 1 ;\,\p\nmid i}
\qint{j}^{
o_\p\left(
\left\lfloor
\frac{m_i(\lm)}{j}
\right\rfloor
\right)}.
\end{align*}
For $(\mmu,\chi)=(\mu^{(1)},\dots,\mu^{(\p-1)},\chi)
\in Q_\p(n)$, define 
$$w_H(\mmu,\chi)=w_H(\mu^{(\p-1)}).$$
Here is an example. 
Suppose that $p=2$ and $\lm = (1^42^345^2)\in P(24)$. 
Then we see that, 
\begin{align*}
w_E(\lm )=&\prod_{i=1,5 }\,\prod_{j=1}^{m_i(\lm)}\,\qnum{j}\\
=&(1)_{[p]}(2)_{[p]}(3)_{[p]}(4)_{[p]}\times (1)_{[p]}(2)_{[p]}\\
=&(2)_{[p]}(4)_{[p]}\times (2)_{[p]}\\
=&[p]_1[p]_1[p]_2\times [p]_1\\
=&[p]_1^3[p]_2, \\
w_H(\lm )=&\prod_{j\geq1} \prod_{i=1,5}\qint{j}^{
o_\p\left(\left \lfloor \frac{m_i(\lm)}{j}\right \rfloor \right)}\\
=&[p]_1^{o_\p\left(\left \lfloor \frac{4}{1}\right \rfloor \right) + o_\p\left(\left \lfloor \frac{2}{1}\right \rfloor \right)}
[p]_2^{o_\p\left(\left \lfloor \frac{4}{2}\right \rfloor \right) + o_\p\left(\left \lfloor \frac{2}{2}\right \rfloor \right)}\\
&\times[p]_3^{o_\p\left(\left \lfloor \frac{4}{3}\right \rfloor \right) + o_\p\left(\left \lfloor \frac{2}{3}\right \rfloor \right)}
[p]_4^{o_\p\left(\left \lfloor \frac{4}{4}\right \rfloor \right) + o_\p\left(\left \lfloor \frac{2}{4}\right \rfloor \right)}
\cdots\\
=&[p]_1^{3+2}[p]_2^{2+1}[p]_3^{1+0}[p]_4^{1+0}[p]_5^{0+0}\cdots\\
=&[p]_1^{5}[p]_2^{3}[p]_3[p]_4. 
\end{align*}
Remark that, when $\p$ is prime, 
it is classically known that $w_E(\lm)|_{q=1}=
\prod_{i\geq 1}(m_i(\lm)!)_\p$ 
and
that $\{ w_E(\lm)|_{q=1}\mid \lm\in P_{(\p)}(n)\}$
 gives the set of the elementary divisors of the 
Cartan matrix of the symmetric group $\mathfrak{S}_n$ in characteristic $\p$.
Remark also that $w_H(\mmu,\chi)|_{q=1}$ gives an expression of block elementary
divisors of $C_n$ (\cite{h}).
\begin{theorem}\label{th;ed}
Let $\p\in\Z_{\geq2}$. Then, as multisets,
\begin{align*}
&\left\{
w_E(\lm)\mid \lm\in P_{(\p)}(n)
\right\}
=
\left\{
w_H(\mmu,\chi)\mid
(\mmu,\chi)\in Q_\p(n)
\right\}.
\end{align*}
\end{theorem}
\noindent{\it Proof.}
\noindent
We write $r=\floor{\frac{n}{\p}}$ throughout this proof.
Put 
$\Pless=\bigsqcup_{0\leq d\leq r}P_{(\p)}(d).$
For a partition $\lm=(1^{l_1}2^{l_2}\dots)\in P_{(p)}(n)$,
define an element of $\Pless$ by
$$
\alpha(\lm)=
(1^{\floor{\frac{l_1}{\p}}}
2^{\floor{\frac{l_2}{\p}}}\dots).
$$
For $(\mmu,\chi)=(\mu^{(1)},\dots,\mu^{(\p-1)},\chi)\in 
M_{\p-1}(d)\times \core_{\p}(n-\p d)\subseteq Q_\p(n)$
with $\mu^{(\p-1)}=(1^{m_1}2^{m_2}\dots)$,
define 
$\beta(\mmu,\chi)$ to be the $\p$-class regular
partition obtained from $\mu^{(\p-1)}$ by removing all parts divisible by $p$, namely,
$$\beta(\mmu,\chi)=\mu^{(\p-1)}\pminus
(\p^{m_\p}2\p^{m_{2\p}}\dots)$$
In this way we obtain the maps 
\begin{align*}
&\alpha: P_{(\p)}(n)\to \Pless,\\
&\beta: Q_\p(n) \to \Pless.
\end{align*}
Moreover, we have
\begin{align*}
&w_E(\lm)=w_H(\alpha(\lm))\ (\hbox{by Lemma}\ \ref{lem;1}\ \hbox{and Lemma}\ \ref{lem;2}),
\\
&w_H(\mmu,\chi)=w_H(\beta(\mmu,\chi)).\ 
\end{align*}
Therefore, to prove the theorem,
it is enough to show that
$$\sharp\alpha^{-1}(\nu)=\sharp\beta^{-1}(\nu)$$
for all $\nu\in \Pless$.

Fix $\nu=(1^{n_1}2^{n_2}\dots)$.
Then
$\alpha\inv(\nu)$
 consists of the elements of the form 
$(1^{\p n_1+e_1}2^{\p n_2+e_2}\dots)$ 
with $0\leq e_i\leq \p-1$.
Since $\sum_{i\geq1}i(\p n_i+e_i)=n$, 
it follows that 
$\alpha\inv(\nu)$ 
is in one-to-one correspondence 
with the set  
$$A(n)=\{
(1^{e_1}2^{e_2},\dots)\in 
P_{(\p )}(n-\p |\nu|)
\mid 0\leq e_i\leq \p -1\hbox{ for }i\geq 1\}
$$
Recall that the generating function
of the set of $\p$-class regular partitions is 
$\frac{\phi(x^\p)}{\phi(x)}$ (the formula \eqref{eq;gen_Pp}).
We have
\begin{align}
\label{eq;genA}
\sum_{n\geq0}\sharp A(n)x^n=
x^{\p|\nu|}\prod_{i\geq1;\,p\nmid i}{\sum_{0\leq e_i<p}{x^{ie_i}}}=
x^{\p|\nu|}\frac{\phi(x^\p)^2}{\phi(x)\phi(x^{\p^2})}.
\end{align}
It is easy to see that
any element of the set $\beta\inv(\nu)$ 
is of the form
$$((\mu^{(1)},\dots,\mu^{({\p}-2)},
\nu\pplus {\p}\tau),\chi)$$ 
for some 
$(\mu^{(1)},\dots,\mu^{({\p}-2)})
\in M_{{\p}-2}(d-(|\nu|+{\p}j)),\ 
\tau=(1^{t_1}2^{t_2}\dots)\in P(j)$
and $\chi\in
\core_{\p}(n-d{\p})$.
Therefore  the set $\beta\inv(\nu)$  
is in one-to-one correspondence with
\begin{align*}
B(n)=\bigsqcup_{0\leq d\leq r}
\bigsqcup_{j\geq0}
M_{{\p}-2}(d-(|\nu|+{\p}j))\times P(j)
\times\core_{\p}(n-d{\p}).
\end{align*}
Using the formulas 
\eqref{eq;gen_M}, \eqref{eq;gen_core} and (2.4), 
we have
\begin{equation}
\label{eq;genB}
\sum_{n\geq0}\sharp B(n) x^n=
x^{\p|\nu|}
\frac{1}
{\phi(x^\p)^{\p-2}}\times\frac{1}{\phi(x^{\p^2})}
\times
\frev{\phi(x)}{\phi(x^\p)^\p}.
\end{equation}
This equals $\sum_{n\geq0}\sharp A(n) x^n$ 
\eqref{eq;genA}
and implies $\sharp A(n)=\sharp B(n)$.
Therefore $\sharp \alpha\inv(\nu)=\sharp\beta\inv(\nu)$.
\qed
\begin{corollary}
$$\prod_{\lambda \in P_{(p)}(n)} w_E(\lambda) 
= \prod_{(\mmu, \chi) \in Q_p(n)} w_H(\mmu,\chi).$$
\end{corollary}
\section
{Determinant formulas and the Glaisher correspondence}
\label{sec;comb}
The elementary divisors and determinant of $\cart$ can be 
expressed by the {\it Glaisher correspondence},
which gives a bijection
between the set of $\p$-regular and $\p$-class regular partitions as described below. 
Let $\lambda=(1^{m_1}2^{m_2}\dots)$ be 
a partition of $n$.  
If $m_i \geq p$, 
then 
transform $\lm$ as
$$ (1^{m_1}2^{m_2}\cdots i^{m_i}
\cdots (pi)^{m_{pi}}\cdots)
 \stackrel{\gstep_i}{\longmapsto}
(1^{m_1}2^{m_2}\cdots i^{m_i-p}
\cdots (pi)^{m_{pi}+1}\cdots).$$
Repeat this procedure until all exponents will get to be less than $p$. The resulting partition
$\tilde{\lambda}$ is $p$-regular.  

Here is an example.
Suppose that $p=2$ and $\lambda=(1^9\,3\,5^3) \in P_{(2)}(27)$.  
Then  we have $\tilde{\lambda}= \gstep_5\gstep_4\gstep_2^2\gstep_1^4(\lambda)$:
\begin{align*}
\lambda&=(1^9\,3\,5^3) 
\stackrel{\gstep_1}{\longmapsto} (1^7\,2\,3\,5^3) 
\stackrel{\gstep_1}{\longmapsto} (1^5\,2^2\,3\,5^3) 
\stackrel{\gstep_1}{\longmapsto} (1^3\,2^3\,3\,5^3)
\stackrel{\gstep_1}{\longmapsto} (1\,2^4\,3\,5^3)\\ 
&\stackrel{\gstep_2}{\longmapsto}  (1\,2^2\,3\,4\,5^3)
\stackrel{\gstep_2}{\longmapsto}  (1\,3\,4^2\,5^3)  
\stackrel{\gstep_4}{\longmapsto}  (1\,3\,5^3\,8) 
\stackrel{\gstep_5}{\longmapsto}  (1\,3\,5\,8\,10) 
= \tilde{\lambda}.
\end{align*}
Now attach weights to the Glaisher correspondence above.
For a $p$-class regular partition $\lambda$, 
let the {\it{Glaisher weight}} $w_G(\lambda)$ 
be 
$
\prod_{i\geq1}\qint{i}^{d_i(\lm)}$,
where $d_i(\lm)$ is the number of occurrences of step $g_i$
in constructing $\lm\mapsto\tilde{\lm}$.
In the previous example, we see that
$w_G(\lambda) = \qint1^4\qint2^2\qint4\qint5$.

For any $\lm\in P_{(p)}(n)$, we have an explicit formula
\begin{align*}
&w_G(\lm)=\prod_{a\geq 1}\prod_{b\geq1}
\qint{a\p^{b-1}}^{\left\lfloor
\frac{m_a(\lm)}{p^{b}} \right\rfloor}.
\end{align*}
Under the specialization $q=1$,
it is not difficult to see that
\begin{equation}
\label{eq;G=E_sp}
w_G(\lm)|_{q=1}=\prod_{i\geq 1}\prod_{b\geq1}
p^{\left\lfloor
\frac{m_i}{p^b} \right\rfloor}
=\prod_{i\geq 1}\,\prod_{j=1}^{m_i}
(j)_p
=w_E(\lm)|_{q=1},
\end{equation}
and that the left hand side is, by definition, equal to  
$p^{ \frac{l(\lambda)-l(\tilde{\lambda})}{p-1}}$,
where $l(\lambda)$ denotes the length of the partition $\lambda$ ([UY]).
It is known that (4.1) gives the elementary divisors of the Cartan
matrix $C_{n}$, and hence the product gives $\det C_{n}$  (e.g. [LLT], [NT]).
%

The  Cartan matrix $\cart$ is also related to the
Gram matrix of the  Shapovalov form for the basic representation 
of the affine Lie algebra $\widehat{\mathfrak{sl}}_{\p}$.
Following Tsuchioka \cite{t}, set
 $$A_j(d)=\sum_{\lm
\in P(d)}
{\frac{m_j(\lm)}{{\p}-1}
\prod_{i\geq1}
\begin{pmatrix}
{\p}-2+m_i(\lm)
\\ m_i(\lm)
\end{pmatrix}},$$
and define
\begin{align*}
&\bd_{\p,n}(d)=\prod_{j\geqq 1}
[\p]_j^{A_j(d)}, \\
&\bd_{\p,n}=\prod_{
0\leq d\leq \floor{\frac{n}{\p}}
}
\bd_{\p,n}(d)^{c_\p(n-\p d) }.
\end{align*}
As shown in \cite{bk1}, 
$\bd_{\p,n}(d)|_{q=1}$ equals the determinant 
of the block of $\p$-weight $d$ of the Cartan matrix $\cart$, and 
it follows that
$\bd_{\p,n}|_{q=1}$
gives the full determinant of $\cart$.
\begin{theorem}\label{th;comb2}
The following equalities hold as polynomials in $q$:
$$
\bd_{\p,n}
=\prod_{\lm\in P_{(\p)}(n)}w_G(\lm)
=\prod_{\lm
\in P_{(\p)}(n)}w_E(\lm).
$$
\end{theorem}
The proof of Theorem \ref{th;comb2} will be
given in the next section.
By comparing the 
exponents
of $[\p]_j$ in
$\prod_{\lm
\in P_{(\p)}(n)}w_G(\lm)$
and 
$\prod_{\lm
\in P_{(\p)}(n)}w_E(\lm)$, 
we obtain the following formula.
\begin{corollary}\label{cor;comb2}
Let $j$ and $k$ be positive integers satisfying $p\nmid j$.
Then 
\begin{align*}
\sum_{\lambda\in P_{(p)}(n)} \left\lfloor\frac{m_{j}(\lm)}{p^{k}}\right\rfloor
=& \sum_{\lambda
\in P_{(p)}(n)} \sum_{i\geq1}
o_\p\left(
\left\lfloor
\frac{m_i(\lm)}{p^{k}j}
\right\rfloor
\right).
\end{align*}
\end{corollary}
It will turn out that the polynomials in Theorem \ref{th;comb2}
coincide with the determinant of the 
{\it graded Cartan matrix},
and
that the following theorem gives an expression
for its block determinant (Theorem \ref{th;blockdet}).
\begin{theorem}\label{th;comb1}
Let $\p\in\Z_{\geq2},\ d\in \Z_{\geq0}$
and $j\in\Z_{\geq1}. $
Then
\begin{align*}
A_j(d)&=
\sum_{
(\mu^{(1)},\dots,\mu^{(\p-1)})\in M_{\p-1}(d)
}
{m_j(\mu^{(\p-1)} )
}\\
&=
\sum_{
(\mu^{(1)},\dots,\mu^{(\p-1)})\in M_{\p-1}(d)}
\sum_{
i\geq 1;\,\p\nmid i}
{\ord_\p 
\left(\left\lfloor
\frac{m_i(\mu^{(\p-1)})}{j}\right\rfloor\right)}.
\end{align*}
Therefore,
$$
\bd_{\p,n}(d)=
\prod_{(\mu^{(1)},\dots,\mu^{(\p-1)})\in M_{\p-1}(d)}w_H(\mu^{(\p-1)}).$$
\end{theorem}
The proof of Theorem \ref{th;comb1}
will be given in Section \ref{sec;pfcomb1}.
Though our proof of Theorem \ref{th;comb2} is 
bijective and the proof of Theorem \ref{th;comb1} is  based on
Theorem \ref{th;comb2},
it is also possible to prove Theorem \ref{th;comb2} and Theorem \ref{th;comb1}
directly and independently by using generating functions. 

The following corollary follows from the second equality in Theorem \ref{th;comb1}.
\begin{corollary}\label{cor;comb1}
Let $j\in\Z_{\geq1}$.
Then
$$\sum_{\lm\in P(n)}m_j(\lm)=\sum_{\lm\in P(n)}\sum_{i\geq1,\; 
\p\nmid i}
\ord_\p \left(\floor{\frac{m_i(\lm)}{j}}
\right).$$
\end{corollary}
\section
{Proof of Theorem \ref{th;comb2}}
\label{sec;glaisher}
To prove the first equality in 
Theorem \ref{th;comb2},
put
\begin{align*}
N_{j,n}&=\sum_{d\geq 0} c_p(n-dp) A_j(d)\\
& =
\sum_{d\geq 0} c_p(n-dp) \left(\sum_{\lambda \in P(d)} \frac{m_{j}(\lm)}{p-1} 
\prod_{i \geq 1} \binom{p-2+m_i(\lm)}{m_i(\lm)} \right).
\end{align*}
Recall that  $c_p(k)=\sharp \core_p(k)$ and
their generating function is given by
$\frev{\phi(x)}{\phi(x^\p)^\p}$.
Note that
$$
\frac{m}{\p-1}\binom{p-2+m}{m} = 
\frac{m}{\p-1}\left(\binom{p-1}{m}\right)
 =\left(\binom{p}{m-1}\right),
$$
where the symbol ``$\left( \left(\,\,\right) \right)$" stands for  the number of combinations 
with repetitions.  It is easy to see that
$$\sum_{m \geq 0} \left(\binom {p-1}{m}\right)x^m = \frac{1}{(1-x)^{p-1}}_.$$
Then it is seen that the generating function of the sequence $\{N_{j,n}\}_{n \geq 0}$ reads
\begin{align*}
\sum_{n\geq 0}N_{j,n}x^n &= 
\frev{\phi(x)}{\phi(x^\p)^\p}
\times \frac{x^{j \p}}{(1-x^{j \p})^\p}\times
\prod_{i \neq j}\frac{1}{(1-x^{i\p})^{\p-1}}\\
&= \frac{x^{j \p}}{1-x^{j \p}}\times \prod_{k\geq 1;\,\p \nmid k}
\frac{1}{1-x^k}_.
\end{align*}
The right-hand side is equal to the generating function
$$\sum_{n \geq 0}\left(\sum_{k \geq 1}\sharp P_{(p)}(n-pk j) \right) x^n.$$
Let $j=ap^{b-1}$, where 
$p$ does not divide $a$, and $b \geq 1$.
Then 
$$\sum_{k \geq 1}\sharp
P_{(p)}(n-pk j) = 
\sum_{\lambda \in P_{(p)}(n)}\left\lfloor\frac{m_{a}}{p^b}\right\rfloor,$$
and this completes the proof of the first equality .

\medskip
Next, we prove the second equality
$$\prod_{\lambda \in P_{(p)}(n)}w_G(\lambda)  = \prod_{\lambda \in P_{(p)}(n)}w_E(\lambda).$$
Our proof is bijective. To this end, we reformulate 
the two weights $w_G$ and $w_E$ as follows.
For $\lambda=(1^{m_1}2^{m_2}\ldots)  \in P_{(p)}(n)$ and $i\geq1;\,\p\nmid i$,
we associate a {\it{diagram}}
$$D_i(\lambda)=\left\{(j,k)\in \mathbb{Z}_{\geq 0} 
\times \mathbb{Z}_{\geq 1}
\;\left|\; 1 \leq 
k \leq \left\lfloor\frac{m_i}{p}\right\rfloor, \right. \,p^j|k\right\}.$$

Here is an example.  Let $p=2$ and $\lambda=(1^9\,3\,5^3) \in P_{(2)}(27)$.
Then we see that 
\def\bx{\framebox(10,10)}
\begin{center}
\setlength{\unitlength}{.05cm}
\begin{picture}(150,60)(-60,-25)
\put(-35,13){$D_1(\lambda)\ =\ $}
\put(0,20){\bx{}}
\put(10,20){\bx{}}
\put(20,20){\bx{}}
\put(30,20){\bx{}}
\put(10,10){\bx{}}
\put(30,10){\bx{}}
\put(30,00){\bx{}\ ,}
\put(-35,-17){$D_5(\lambda)\ =\ $}
\put(0,-20){\bx{}\ ,}
\end{picture}
\end{center}

\noindent
and $D_i(\lambda) = \emptyset$  for other odd $i$.

Put
$$\mathfrak{D}(n) =\mathfrak{D}(n,p) =\left\{
(\lambda; i,j,k)\mid \lambda \in P_{(p)}(n),\, 
\p\nmid i,
\, 
(j,k) \in D_i(\lambda)\right\}.$$
We consider two {\it{tableaux}} $G$ and $E$ on $\mathfrak{D}(n)$, namely
$$G, \, E : \mathfrak{D} (n) \longrightarrow \mathbb{Z}_{\geq 1}.$$
For $c=(\lambda; i, j,k) \in \mathfrak{D}(n)$, define $G(c)=ip^{j}$ and $E(c) = k/p^j$ .

In the previous example, $G$ and $E$ are tabulated on $D_1(\lambda)$, respectively, as
\begin{center}
\setlength{\unitlength}{.05cm}
\begin{picture}(150,40)(-60,-5)
\put(-45,13){$G(D_1(\lambda))\ =\ $}
\put(0,20){\bx{$1$}}
\put(10,20){\bx{$1$}}
\put(20,20){\bx{$1$}}
\put(30,20){\bx{$1$}}
\put(10,10){\bx{$2$}}
\put(30,10){\bx{$2$}}
\put(30,00){\bx{$4$}\ ,}
\end{picture}
\end{center}
\begin{center}
\setlength{\unitlength}{.05cm}
\begin{picture}(150,40)(-60,-5)
\put(-45,13){$E(D_1(\lambda))\ =\ $}
\put(0,20){\bx{$1$}}
\put(10,20){\bx{$2$}}
\put(20,20){\bx{$3$}}
\put(30,20){\bx{$4$}}
\put(10,10){\bx{$1$}}
\put(30,10){\bx{$2$}}
\put(30,00){\bx{$1$}\ .}
\end{picture}
\end{center}
\noindent
Here $G(\lambda; 1,j,k)$ (resp. $E(\lambda; 1,j,k)$) is written in the
$(j,k)$-position.

\medskip

The claim of the theorem is equivalent to
the identity
\begin{equation}\label{eq;p=p}
\prod_{c \in \mathfrak{D}(n)} \qint{G(c)} = \prod_{c \in \mathfrak{D}(n)} \qint{E(c)}. 
\end{equation}
We prove \eqref{eq;p=p} by constructing an involution
$$\theta : \mathfrak{D}(n) \longrightarrow \mathfrak{D}(n)$$
such that $E\circ \theta = G$.
Take an element  $c=(\lambda; i,j,k) \in \mathfrak{D}(n)$.  By definition,  $pk \leq m_i$ and
$p^j|k$.  Therefore there exists a $p$-class regular partition $\mu$ such that 
$\lambda= \mu + (i^{pk})$,
where $``+"$ denotes the concatenation (union) of two Young diagrams.  
Let $k$ be decomposed as
$k=i'p^{j+j'}$, where $i'=(k)_{p'}$ is the $p'$-part of $k$, and $p^{j+j'}=(k)_p$ is the
$p$-part of $k$.  
Hence we can write 
$$c=(\mu + (i^{pi'p^{j+j'}}); i,j,i'p^{j+j'}).$$
Define $\theta$ by 
$$\theta(c) = (\mu + (i'^{pip^{j+j'}}); i',j', ip^{j+j'}).$$
Namely, $\theta$ interchanges $i$ with $i'$, and $j$ with $j'$.  It can be seen that 
$\alpha$ is a map 
from $\mathfrak{D}(n)$ to itself, and is an involution.  It is also easy to see that, for
$c=(\lambda; i,j,k)$, 
$$E(\theta(c)) = ip^j = G(c).$$
This proves the formula.
\section{Proof of Theorem \ref{th;comb1}
}\label{sec;pfcomb1}
\noindent{\it Proof of the first equality.}
For each $j$, we have
\begin{align*}
\sum_{d\geq0} \sum_{\mmu\in M_{{\p}-1}(d)} 
m_j(\mu^{({\p}-1)})x^d
&=\sum_{d\geq0}\sum_{k=0}^d\sum_{(\mlm,\mu)\in 
M_{{\p}-2}(d-k)\times P(k)}
m_j(\mu) x^d\\
&=\frac{1}{\phi(x)^{{\p}-2}}
\times\frac{1}{\phi(x)}\times \frev{(1-x^j)}{{x^j}}
=\frac{x^j}{(1-x^j){\phi(x)^{{\p}-1}}},
\end{align*}
\begin{align*}
\sum_{d\geq0}\sum_{\lm\in P(d)}
\frac{m_j}{{\p}-1}
\prod_{i\geq1}
\left(
\begin{pmatrix}
{\p}-1\\ m_i
\end{pmatrix}
\right)
x^{d}
=&
\frac{x^j}{(1-x^j)^{\p}}\times
\prod_{i\geq1,i\neq j}\frac{1}{(1-x^i)^{{\p}-1}}\\
=&\frev{(1-x^j)\phi(x)^{{\p}-1}}{x^j}.
\end{align*}
Hence the formula follows.
\qed

\bigskip\noindent
{\it Proof of the second equality.}

\smallskip

By Theorem \ref{th;ed} and Theorem \ref{th;comb2}, we have 
\begin{align*}
\prod_{(\mlm,\chi)\in Q_\p(n)}
w_H(\lm^{(\p-1)})
&=\prod_{\lm\in P_{(\p)}(n)}w_E(\lm)
=\bd_{\p,n}
\end{align*}
By comparing the exponent of $\qint{j}$,
we have
\begin{align*}
&\sum_{0\leq d\leq \lfloor \frac{n}{\p} \rfloor }
\coren_p(n-\p d)\delta_d=0,
\end{align*}
for each $n$,
where
$$\delta_d=\sum_{\underline{\lm}\in M_{\p-1}(d)}
\left(
A_j(d)-
\sum_{i\geq1;\,\p\nmid i}
\ord_\p \floor{\frac{m_i(\lm^{(\p-1)})}{j}}
\right)=0.$$
 By letting $n=\p d'$, we have
$$
\delta_{d'}+\sum_{0\leq d<d'}\coren_p(n-\p d)
\delta_d=0
$$
as 
$\coren_p(0)=1$.
By induction on $d'$, we have $\delta_{d'}=0$, namely,
\begin{align*}
&A_j(d)=
\sum_{\underline\lm\in M_{\p-1}(d)}
\sum_{i\geq1;\,\p\nmid i}
\ord_\p \floor{\frac{m_i(\lm^{(\p)})}{j}}.
\end{align*}
\qed
\section{Block version for $p=2$}
\label{sec;p=2}
When $p=2$, the Glaisher correspondence turns out to be a bijection from the set of the
odd partitions to the set of the strict partitions.  In this case we can refine Theorem 4.1 to 
block-wise products.  We write
$SP(m) = P^{(2)}(m),\,\,OP(m) = P_{(2)}(m)$ and $OSP(m) = OP(m) \cap SP(m)$.
We use the following diagram representing a strict partition, 
which is called the {\it{4-bar abacus}} in \cite{bo} and the {\it{H-abacus}} in \cite{uy}. 
For example,
the H-abacus of $\lambda = (2,3,7,9)$ is shown below.

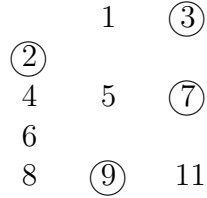
\begin{center}
\begin{picture}(80,80)
\put (40,65){1}
\put (70,65){3}
\put (10,50){2}
\put (10,35){4}
\put (40,35){5}
\put (70,35){7}
\put (10,20){6}
\put (10,5){8}
\put (40,5){9}
\put (68,5){11}
\put (12.5,53){\circle{13}}
\put (72.5,68){\circle{13}}
\put (72.5,38){\circle{13}}
\put (42.5,8){\circle{13}}
\end{picture}
\end{center}

\noindent
Namely, for a strict partition we put a set of beads on the assigned
positions.  
Two beads do not occupy the same position. 
>From the H-abacus of
the given strict partition $\lambda$, we obtain the
{\it{H-core}} $\lambda^H$ by moving and removing the beads as follows:

\vspace{5mm}

(1)  Move a bead one position up along the leftmost runner.

(2)  Remove a bead at the position 2.

(3)  Move a bead one position up along the runner of 1 or of 3.

(4)  Remove the two beads at the positions 1 and 3, simultaneously.
\vspace{5mm}

\noindent
The H-cores are thus characterized by the {\it{stalemates}}, which
constitute the set
$$HC = \left\{\emptyset,\, (1,5, \ldots, 4m-3, 4m+1),\, (3,7,\ldots, 4m-1,4m+3)\mid
\, m \geq 0\right\}.$$
For example, the H-core of the above $\lambda=(2,3,7,9)$ is
$\lambda^H = (3)$.
Notice that the number of nodes in every H-core is a triangular number, $m(m+1)/2$, and conversely, for any triangular number $r$, there is a unique H-core with $r$ nodes. Thus there is a unique bijection between $HC$ and the set of 2-cores 
$$\left\{\Delta_0 = \emptyset,\, \Delta_m = (1,2,\ldots, m) \mid m\geq 1\right\}$$
that preserves the number of nodes. In fact, the bijection is obtained by applying {\it{unfolding}}, which is defined as taking the hook lengths
of the main diagonal in the Young diagram. Namely, we have
$$HC = \{\Delta_m^{u}\mid m \geq 0\}, $$
where  $\lambda^{u}$
stands for the unfolding of $\lambda$.  For example,
$\Delta_4^{u} = (3,7)$.

We need the {\it{H-quotient}} $\lambda^H[1]$ for a strict partition $\lambda$.
Draw the H-abacus of the strict partition $\lambda$, and read out a 0-1 sequence
as follows. First look at the runner of 3 starting from the bottom.  If the number is circled, then
attach 0, and attach 1 otherwise.  In the above example, $\lambda = (2,3,7,9)$, 
this 0-1 sequence is ...100.  Next look at the runner of 1 starting at the top.  If the number
is circled, then attach 1 and attach 0 otherwise. In the above example,
the 0-1 sequence is 0010... .  Concatenate two 0-1 sequences to get
a two-side infinite 0-1 sequence. In the above example we have
...1000010...  .  From this {\it{Maya diagram}} we define a partition
$\lambda^H[1]$ by counting
0's on the left of each 1.  In the example, we have the partition $\lambda^H[1] = (4)$.
It should be noticed that, for every fixed H-core $\lambda^H$, the map
$\lambda \mapsto \lambda^H[1]$ is a bijection from $OSP(4d+|\lambda^H|)$ to
$P(d)$. (See for example \cite{ol}.)

\begin{theorem}
For any non-negative integer $d$,
$$\prod_{\lambda}w_G(\lambda)  = \prod_{\lambda}w_E(\lambda),$$
where the products of the both sides run over all odd partitions $\lambda$ of $2d$ such that
$\tilde{\lambda}^H = \emptyset$.
Moreover, they equal the block determinant of the graded Cartan matrix of $2$-weight $d$.
\end{theorem}
\noindent{\it Proof.}
This first claim is easily verified by noticing that 
the involution 
$\theta$ does not change the H-core
of $\tilde{\lambda}$.  In fact, according to the decomposition
$\lambda = \mu + (i^{pk})$ in the proof of 
Theorem \ref{th;comb2} 
$\tilde{\lambda}^H = \tilde{\mu}^H$.

As for the second half we show 
\begin{align*}
&\left\{
w_E(\lm)\mid \lm\in OP(2d), \tilde{\lambda}^H = \emptyset
\right\}
=
\left\{
w_H(\mu,\emptyset)\mid
(\mu,\emptyset)\in Q_2(d)
\right\}.
\end{align*}
as multisets.

This can be shown in a similar way to that of 
Theorem \ref{th;ed}.
Put $r = \lfloor \frac{n}{2} \rfloor$ and recall the two maps
\begin{align*}
&\alpha: OP(n)\to OP(\leq r)\\
&\beta: Q_2(n) \to OP(\leq r). 
\end{align*}
Restrict these maps to the subsets
$A'(d) = \{\lambda \in OP(2d)\mid \tilde{\lambda}^H = \emptyset\}$ and
$B'(d) = \{\mu \in P(d)\mid (\mu, \emptyset) \in Q_2(d)\} $, respectively, and keep the same notation
\begin{align*}
&\alpha : A'(d) \to OP(\leq d)\\
&\beta: B'(d) \to OP(\leq d)
\end{align*}
Here we identify the pair $(\mu, \emptyset) \in Q_2(d)$ with a single partition $\mu \in P(d)$,
and write $w_H(\mu)$ in place of $w_H(\mu, \emptyset)$.

It is easily seen that $w_E(\lambda) = w_H(\mu)$ if $\alpha(\lambda) = \beta(\mu)$.
Therefore it is enough to show that 
$$\sharp\alpha^{-1}(\nu)=\sharp\beta^{-1}(\nu).$$
for all $\nu\in OP(\leq d)$.

Fix $\nu=(1^{n_1}2^{n_2}\dots)$.
Then
$\alpha^{-1}(\nu)$
 consists of the elements of the form 
$(1^{2 n_1+e_1}2^{2 n_2+e_2}\dots)$ 
with $e_i = 0$ or $1$.
Since $\sum_{i\geq1}i(2n_i+e_i)=2d$, 
it follows that 
$\alpha^{-1}(\nu)$ 
is in one-to-one correspondence 
with the set
$OSP(2d-2|\nu|)$.

On the other hand, 
any element of the set $\beta^{-1}(\nu)$ 
is of the form
$\nu + 2\tau$
for some 
$\tau=(1^{t_1}2^{t_2}\dots)\in P(\frac{d-|\nu|}{2})$.
It is already mentioned that the two sets
$OSP(2d-2|\nu|)$ and $ P(\frac{d-|\nu|}{2})$ have the same 
cardinality.
\qed
\section{Graded Cartan matrices}
\label{sec;qcartan}
%
In 1996, Lascoux, Leclerc and Thibon \cite{llt} 
presented an algorithm for computing the global
crystal basis for the basic representation $L(\Lambda_0)$
of $U_q({\widehat{\mathfrak{sl}}_p})$,
where $q$ is an indeterminate.
The basic representation of $U_q({\widehat{\mathfrak{sl}}_p})$ is realized as the highest irreducible component of the Fock space
$$\mathfrak{F} = \bigoplus_{\lambda \in P} \mathbb{Q}(q)\lambda$$
where $P$ denotes the set of all partitions.  
Let $d_{\lambda \mu}(q)$ be determined by
$$\mathcal{G}(\mu) = \sum_{\lambda\in P}d_{\lambda \mu}(q)\lambda,$$
where $\{\mathcal{G}(\mu)\mid \mu \in P^{(p)}=\bigsqcup_{n\geq0}P^{(p)}(n)\}$ 
is Kashiwara's {\it {lower global crystal basis}} 
for the basic representation of $U_q({\widehat{\mathfrak{sl}}_p})$.  
Define the matrices
\begin{align*}D_n(q) &= 
(d_{\lambda \mu}(q))_{\lm\in P(n),\, \mu\in P^{(p)}(n)},\\ 
C_n(q)&={}^t D_n(q)D_n(q).
\end{align*}
Lascoux, Leclerc and Thibon \cite{llt} conjectured that
$D_n(1) = (d_{\lambda \mu}(1))$ is the decomposition matrix of 
the Iwahori-Hecke algebra $H_n(\zeta)$ of type $A_{n-1}$
with $\zeta$ being a primitive $p$-th root of unity, and
this conjecture was proved by Ariki \cite{a}.

A representation-theoretic background of
$D_n(q)$ is given by
the {\it Khovanov-Lauda-Rouquier algebra}
$\KLR$ associated with the symmetric group,
introduced independently 
by Khovanov and Lauda \cite{kl1,kl2}
and Rouquier \cite{rouq}.
The algebra $\KLR$ is a $\Z$-graded algebra
and is 
isomorphic to 
the Iwahori-Hecke algebra $H_n(\zeta)$
as a non-graded algebra. 
It is shown by 
Brundan and Kleshchev \cite{bk2} 
that 
$D_n(q)$ is nothing but
the ``graded decomposition matrix" for $\KLR$.
As a consequence, it follows that
$C_n(q)$ 
is the corresponding graded Cartan matrix
(see also \cite[Theorem 2.17]{hm}). 
More precisely, 
$$
C_n(q)=\left(
\sum_{k\in\Z}
[\widetilde{P}(\lm):
\widetilde{D}(\mu)\langle k \rangle]q^k
\right)_{\lm,\mu\in P^{(p)}(n)}$$
Here, $\widetilde{D}(\mu)\langle k\rangle$ 
is the graded simple module of $\KLR$
which corresponds to $D(\mu)$
(the simple $H_n(\zeta)$-module of label $\mu$)
 by the grade 
forgetting functor,
$\langle k \rangle$ indicates
the grading shift, and
$\widetilde{P}(\lm)$ denotes the projective cover
of $\widetilde{D}(\lm)\langle 0\rangle$.
 We remark here that 
$C_n(0) = E$, the unit matrix 
\cite[Theorem 6.8]{llt}.

\def\Uq{U_q(\widehat{\mathfrak{sl}}_q(p))}
Using the description of the Fock space in 
\cite[Theorem 3]{cj},
Tsuchioka \cite{t} obtained a formula
for the block determinants of the Gram matrix of 
the Shapovalov form for the basic $U_q(\widehat{\mathfrak{sl}}_p)$-module
by a similar argument in \cite{bk1}.

As is well-known, the blocks are labeled by the $p$-cores.
A partition $\lambda$ is said to have weight $d$ if $d$ successive removals of 
$p$-hooks from $\lambda$ achieve a $p$-core.  On the other hand, 
the block of weight $d$ of the Cartan matrix corresponds to the 
Gram matrix of the Shapovalov form for the weight space of weight
$w\Lambda_0 - d\delta$, where $w$ is an element of the Weyl group
and $\delta$ is the fundamental imaginary root 
of $\widehat{\mathfrak{sl}}_p$
 (see for example \cite{bh}).

We remark that 
the Shapovalov determinant for $L(\Lambda_0)$
is determined up to 
the even powers of $q$, 
and we  normalize the determinant so that 
it equals 1 at $q=0$.  
%

\begin{theorem}[Tsuchioka]
\label{th;blockdet}
Let $d$ be 
the weight of the block. 
Then the corresponding 
block determinant of the graded Cartan matrix $C_n(q)$ is given by
$$\prod_{j\geq 1} \qint{j}^{A_j(d)},$$
where
$$A_j(d) = \sum_{\lambda \in P(d)} \frac{m_{j}(\lm)}{\p-1} 
 \prod_{i \geq 1}\binom{p-2+m_i(\lm)}{m_i(\lm)}$$
 as before.
\end{theorem}
%
Tsuchioka's formula gives a natural graded 
analogue of the 
block determinant formula for the Cartan matrix 
$C_n$ by Brundan and Kleschchev \cite{bk1}.
As an analogue of the expression for the
elementary divisors for $C_n$, 
we conjecture the following:\\

\begin{conjecture}\label{con}
Let $\p$ be a prime.

\smallskip\noindent
$\mathrm{(i)}$ The elementary divisors of the block of
$\p$-weight $d$ of the graded Cartan matrix 
over the ring $\mathbb{Q}[q,q^{-1}]$ coincide with the 
elementary divisors of the diagoal matrix with entries 
$$
\left\{
w_H(\mu^{(\p-1)})\mid
\mmu=(\mu^{(1)},\dots,\mu^{(\p-1)})\in M_{\p-1}(d)
\right\}.$$

\smallskip\noindent
$\mathrm{(ii)}$ 
The elementary divisors of the graded Cartan matrix $C_{n}(q)$ 
over the ring $\mathbb{Q}[q,q^{-1}]$ 
coincide with the elementary divisors of the diagonal matrix with entries 
$$\left\{
w_E(\lm)\mid \lm\in P_{(\p)}(n)
\right\},$$
and also that of the diagonal matrix with entries 
$$\left\{
w_G(\lm)\mid \lm\in P_{(\p)}(n)
\right\}.$$
\end{conjecture}


\bigskip
\noindent
{\bf Acknowledgments.}
We thank Shunsuke Tsuchioka (IPMU) 
for explaining his new results.
We also thank Katsuhiro Uno (Osaka Kyoiku University) 
for valuable discussions and suggestions.

\end{document}